\documentclass[11pt]{amsart}
\usepackage{times}
\usepackage{amsfonts}
\usepackage{amsthm}
\usepackage{amsmath}
\usepackage{epic}
\usepackage{times}
\usepackage{verbatim}
\usepackage{a4wide}
\newcommand{\Orbit}{\mathcal{O}}
\newcommand{\Closure}{\overline{\Orbit}}

\def\cA{\mathcal{A}}

\def\empt{\varepsilon}
\def\NN{\mathbb{N}}

\def\bs{\boldsymbol{s}}
\def\bu{\boldsymbol{u}}
\def\bv{\boldsymbol{v}}

\def\bbf{\boldsymbol{f}}

\def\bx{\boldsymbol{x}}
\def\by{\boldsymbol{y}}
\def\bt{\boldsymbol{t}}
\def\bc{\boldsymbol{c}}

\def\cP{\mathcal{P}}

\def\cP{\mathcal{P}}
\def\cL{\mathcal{L}}
\def\cS{\mathcal{S}}

\theoremstyle{plain}
\newtheorem{theorem}{Theorem} 
\newtheorem{lemma}[theorem]{Lemma}
\newtheorem{corollary}[theorem]{Corollary}
\newtheorem{proposition}[theorem]{Proposition}

\theoremstyle{definition}
\newtheorem{definition}[theorem]{Definition}
\newtheorem{example}[theorem]{Example}

\theoremstyle{remark}
\newtheorem{remark}[theorem]{Remark}
\newtheorem*{note}{Note}

\begin{document}

\title{Distribution modulo 1 and the lexicographic world}

\author{Jean-Paul Allouche}

\address{CNRS, LRI, UMR 8623, Universit\'e Paris-Sud, B\^atiment 490, 
F-91405 Orsay Cedex, FRANCE}

\email{allouche@lri.fr}

\author{Amy Glen}

\address{Department of Mathematics and Statistics, School of Chemical and Mathematical 
Sciences, Murdoch University, Perth, WA 6150, AUSTRALIA}

\curraddr{The Mathematics Institute, Reykjav\'ik University, Kringlan 1, 
IS-103 Reykjav\'ik, ICELAND}

\email{amy.glen@gmail.com}



\dedicatory{In honour of Paulo Ribenboim on the occasion of his 80th birthday}

\begin{abstract}
We give a complete description of the minimal intervals 
containing all fractional parts $\{\xi 2^n\}$, for some positive real number $\xi$,
and for all $n \geq 0$.
\end{abstract}

\maketitle

\section{Introduction}

In the paper \cite{kM68anun} Mahler defined the set of $Z$-numbers by
$$
\left\{
\xi \in {\mathbb R}, \ \xi > 0, \ \forall n \geq 0, \
0 \leq \left\{ \xi \left(\frac{3}{2}\right)^n\right\} < \frac{1}{2}
\right\}
$$
where $\{z\}$ is the fractional part of the real number $z$.
Mahler proved that this set is at most countable. It is still an open problem to
prove that this set is actually empty. More generally, given a real number $\alpha > 1$
and an interval $(x, y) \subset (0,1)$ one can ask whether there exists $\xi > 0$
such that, for all $n \geq 0$, we have $x \leq \{\xi \alpha^n\} < y$ (or the variant
$x \leq \{\xi \alpha^n\} \leq y$).
Flatto, Lagarias, and Pollington \cite[Theorem~1.4]{FLP} proved that, if $\alpha = p/q$
with $p, q$ coprime integers and $p > q \geq 2$, then any interval $(x, y)$ such that
for some $\xi > 0$, one has that $\{\xi (p/q)^n\} \in (x, y)$ for all $n \geq 0$,
must satisfy $y-x \geq 1/p$. Recently Bugeaud and Dubickas \cite{yBaD05frac} 
characterized irrational numbers $\xi$ such that for a fixed integer $b \geq 2$
all the fractional parts $\{\xi b^n\}$ belong to a closed interval of length $1/b$.
Before stating their theorem we need a definition.

\begin{definition} \label{D:intro}
Given two real numbers $\alpha$ and $\rho$ with $\alpha \geq 0$, we denote by 
$\bs_{\alpha,\rho} := (s_{\alpha,\rho}(n))_{n \geq 0}$ and  
$\bs'_{\alpha,\rho} := (s'_{\alpha,\rho}(n))_{n \geq 0}$ the sequences defined by
$$
s_{\alpha,\rho}(n) = \lfloor(n+1)\alpha + \rho\rfloor -
                        \lfloor n\alpha + \rho\rfloor \quad \mbox{and} \quad 
s'_{\alpha,\rho}(n) = \lceil(n+1)\alpha + \rho\rceil -
                        \lceil n\alpha + \rho\rceil \quad \mbox{for $n\geq 0$},
$$
where $\lfloor x \rfloor$ denotes the greatest integer $\leq x$ and $\lceil x \rceil$ denotes 
the least integer $\geq x$. The sequences $\bs_{\alpha,\rho}$ and $\bs'_{\alpha,\rho}$ are 
called {\em Sturmian sequences} if $\alpha$ is irrational, and {\em periodic balanced sequences} 
if $\alpha$ is rational. Furthermore, if $\alpha = \rho$, these sequences are called 
{\em characteristic} Sturmian sequences or {\em characteristic} periodic
balanced sequences, according to whether or not $\alpha$ is irrational.
\end{definition} 

We observe that, if $\alpha$ is not an integer, then for all $n \geq 0$, 
\[
\lfloor \alpha \rfloor \leq s_{\alpha,\rho}(n) \leq \lfloor \alpha \rfloor + 1 
\quad \mbox{and} \quad \lceil \alpha \rceil - 1 \leq s'_{\alpha,\rho}(n) \leq 
\lceil \alpha \rceil,
\]
where $\lceil \alpha \rceil -1 = \lfloor \alpha \rfloor$ and 
$\lceil \alpha \rceil = \lfloor \alpha \rfloor + 1$. On the other hand, 
if $\alpha$ is an integer, then $\bs_{\alpha,\rho} = \bs'_{\alpha,\rho}$ 
and $\alpha \leq s_{\alpha, \rho}(n) \leq \alpha + 1$ for all $n \geq 0$. 
Accordingly, the sequences $\bs_{\alpha,\rho}$ and $\bs'_{\alpha,\rho}$ 
take their values in the ``alphabet'' $\{k, k+1\}$ where $k = \lfloor \alpha \rfloor$.  
The classical definition of Sturmian sequences with values in $\{0,1\}$ is thus 
obtained by subtracting $\lfloor \alpha \rfloor$ from each of the terms in the 
sequences $\bs_{\alpha,\rho}$ and $\bs'_{\alpha,\rho}$. Alternatively, one may 
restrict $\alpha$ to the interval $(0,1)$. Hereafter, if the alphabet is not mentioned, 
it is understood that the sequences are over $\{0, 1\}$. We may also assume that 
$\rho \in [0,1)$ or $\rho \in (0,1]$ since $\bs_{\alpha,\rho} = \bs_{\alpha,\rho'}$ 
and $\bs'_{\alpha,\rho} = \bs'_{\alpha,\rho'}$ for any two real numbers $\rho$, $\rho'$ 
such that $\rho - \rho'$ is an integer.

\begin{example} \label{Ex:Fibonacci} Taking $\alpha = \rho = (3 - \sqrt{5})/2$, we get 
the well-known (binary) {\em Fibonacci sequence} $0100101001001010010100100101\cdots$. 
\end{example}

\begin{remark} \label{R:slope} Note that if $\alpha$ is irrational, then the (Sturmian) 
sequences $\bs_{\alpha,\rho}$ and $\bs'_{\alpha,\rho}$ are {\em aperiodic} (i.e., not 
eventually periodic), whereas if $\alpha$ is rational, the sequences $\bs_{\alpha,\rho}$ 
and $\bs'_{\alpha,\rho}$ are (purely) periodic. (See for instance  \cite[Lemma~2.14]{mL02alge}.) 
This justifies the use of ``periodic'' in the name of such sequences in the rational case. 
The reason for being called ``balanced'' is explained in Section~\ref{SS:balanced}.
\end{remark}

Let $T$ denote the {\em shift map} on sequences, defined as follows: if
$\bs := (s_n)_{n \geq 0}$, then $T(\bs) = T((s_n)_{n \geq 0}) := (s_{n+1})_{n \geq 0}$. 
The main result in \cite{yBaD05frac} reads as follows.

\begin{theorem}[Bugeaud-Dubickas]\label{bugdub}
Let $b \geq 2$ be an integer and let $\xi$ be an irrational number. Then the numbers
$\{\xi b^n\}$ cannot all lie in an interval of length $< 1/b$. Furthermore there
exists a closed interval $I$ of length $1/b$ containing the numbers $\{\xi b^n\}$
for all $n \geq 0$ if and only if the sequence of base $b$ digits of the fractional
part of $\xi$ is a Sturmian sequence $\bs$ on the alphabet $\{k, k+1\}$ for some
$k \in \{0, 1, \ldots, b-2\}$. If this is the case, then $\xi$ is transcendental, 
and the interval $I$ is semi-open. It is open unless there exists an integer 
$j \geq 1$ such that $T^j(\bs)$ is a characteristic Sturmian sequence on the 
alphabet $\{k, k+1\}$.
\end{theorem}

The purpose of this paper is to give a complete description of the minimal intervals 
containing all fractional parts $\{\xi 2^n\}$ for some positive real number $\xi$,
and for all $n \geq 0$.  More precisely, inspired by the definition of the
{\em lexicographic world} (see Section~\ref{lex}), let us define a function $F$ on $[0, 1]$ as follows.

\begin{definition} For all $x \in [0,1]$, let 
$S_{x} := \{\xi \in \mathbb{R}, \ \xi > 0, \ \forall n \geq 0, \ x \leq \{\xi2^n\} < 1\}$ 
and let $F: [0,1] \rightarrow [0,1]$ be the function defined by:
\[
F(x) = \begin{cases}
             \inf\{y \in [0, 1), \ \exists \xi > 0, \ \forall n \geq 0, \ x \leq \{\xi 2^n\} \leq y\} 
                     &\mbox{if $S_{x} \ne \emptyset$}, \\
            1 &\mbox{if $S_{x} = \emptyset$}.
             \end{cases}
\]
\end{definition}

\begin{remark} From Bugeaud-Dubickas'  result for $b=2$, we deduce the following two facts.
\begin{itemize}
\item For $x \in [\frac{1}{2},1]$, there does not exist an {\em irrational} number $\xi > 0$ 
such that $x \leq \{\xi 2^n\} < 1$ for all $n \geq 0$. Nor does there exist a {\em rational} 
number $\xi > 0$ such that $x \leq \{\xi2^n\} < 1$ for all $n \geq 0$. (This can be seen by 
considering, for instance, the base $2$ expansion of the fractional part of $\xi$ for any 
rational number $\xi \geq x$.) Hence, $F(x) = 1$ for  all $x \in [\frac{1}{2},1]$. 
\item If $\xi > 0$ is an irrational real number, then there exists a real number 
$x \in [0,\frac{1}{2})$ such that all the fractional parts $\{\xi 2^n\}$ belong to the interval 
$[x, x + \frac{1}{2}]$ if and only if the base $2$ expansion of the fractional part of $\xi$ is
a Sturmian sequence. Furthermore, for any such $x$, one has 
$F(x) = x + \frac{1}{2}$.
\end{itemize}
Note that it follows from our main theorem (see Theorem~\ref{main} later) that 
$0 \leq F(x) < 1$ for $x \in [0, \frac{1}{2})$.
\end{remark}

Before stating our main theorem, let us note that the sequences $\bs_{\alpha,\rho}$ 
and $\bs'_{\alpha,\rho}$ (given in Definition~\ref{D:intro}) are said to have {\em slope} 
$\alpha$ and {\em intercept} $\rho$, in view of their geometric realization as approximations 
to the line $y = \alpha x + \rho$ (called {\em lower} and {\em upper mechanical words} in 
\cite[Chapter 2]{mL02alge}). From now on, we will assume  that $\alpha$ and $\rho$ are in 
the interval $[0,1]$, in which case the sequences $\bs_{\alpha,\rho}$ and $\bs'_{\alpha,\rho}$ 
take their values in $\{0,1\}$. If $\alpha$ is irrational, then we have $\bs_{\alpha,\alpha} 
= \bs'_{\alpha,\alpha}$, denoted by $\bc_\alpha$. We also have $\bs_{0,0} = \bs'_{0,0} 
= 0^\infty$ and $\bs_{1,1} = \bs'_{1,1} = 1^\infty$, denoted by $\bc_0$ and $\bc_1$, 
respectively. In these cases, the sequence $\bc_\alpha$ is the unique characteristic Sturmian 
sequence of slope $\alpha$ in $\{0,1\}^\NN$. On the other hand, if $\alpha \in (0,1)$ is 
rational, then the characteristic periodic balanced sequences of slope $\alpha$, namely 
$\bs_{\alpha,\alpha}$ and $\bs'_{\alpha,\alpha}$, are distinct sequences in $\{0,1\}^\NN$ 
containing both $0$'s and $1$'s. More precisely, let us suppose that $\alpha = p/q \in (0,1)$ 
with $\gcd(p,q) = 1$. Then by considering the prefix of length $q$ of each of the sequences 
$\bs_{p/q,0}$ and $\bs'_{p/q,0}$, we find that there exists a unique word $w_{p,q}$ of length 
$q - 2$ such that 
\[
\bs_{p/q,0} = (0w_{p,q}1)^\infty \quad \mbox{and} \quad  \bs'_{p/q,0} = (1w_{p,q}0)^\infty
\]
where $v^\infty$ denotes the periodic sequence $vvvv\cdots$ for a given word $v$. (See for 
instance \cite[pg. 59]{mL02alge}.) Hence, the two characteristic periodic balanced sequences 
of slope $p/q$ in $\{0,1\}^\NN$ are given by 
\[
\bs_{p/q,p/q} = T(\bs_{p/q,0}) = (w_{p,q}10)^\infty \quad \mbox{and} 
\quad \bs'_{p/q,p/q} = T(\bs'_{p/q,0}) = (w_{p,q}01)^\infty.
\]
The words $w_{p,q}$ are often referred to as {\em central words} in the literature; they hold 
a special place in the rich theory of Sturmian sequences (see, e.g., \cite[Chapter 2]{mL02alge}).
For instance, it follows from the work of de~Luca and Mignosi~\cite{aD97stur, aDfM94some} that 
central words coincide with the palindromic prefixes of characteristic Sturmian sequences (see 
Section~\ref{SS:periodic}).

Given a sequence $\bs \in \{0,1\}^\NN$, let $r(\bs)$ denote the real number whose sequence of 
base $2$ digits is given by $\bs$. Our main number-theoretical result reads as follows.

\begin{theorem}\label{main}
Let $x$ be a real number in $[0, 1]$.
\begin{itemize}
\item[(i)] If $x \geq \frac{1}{2}$, then $F(x) = 1$.
\item[(ii)] If $x = 0$, then $F(x) = 0$.
\item[(iii)] If $x \in (0,\frac{1}{2})$ and if the base $2$ expansion of $2x$ 
is given by a characteristic Sturmian sequence, then $F(x) = x + \frac{1}{2}$.
Furthermore, $F(x)$ is the unique real number in $[0,1]$ that has a Sturmian base $2$ 
expansion and satisfies $x \leq \{F(x) 2^k\} \leq F(x)$ for all $k \geq 0$.
\item[(iv)] If $x \in (0, \frac{1}{2})$ and if the base $2$ expansion of $2x$ is given by 
a characteristic periodic balanced sequence of slope $p/q \in (0,1)$ with $\gcd(p,q) = 1$, 
then $F(x)$ is the rational number whose base $2$ expansion is given by the periodic balanced 
sequence $\bs'_{p/q,0} = (1w_{p,q}0)^\infty$, in which case $F(x)~\leq~x~+~\frac{1}{2}$.
\item[(v)] In all other cases, $F(x)$ can be explicitly computed: it is equal to the rational 
number whose base $2$ expansion is given by a (unique) periodic balanced sequence 
$\bs'_{p/q,0} = (1w_{p,q}0)^\infty$ where $p$, $q$ are coprime integers with $0< p < q$ such 
that $r((w_{p,q}01)^\infty) < 2x < r((w_{p,q}10)^\infty)$. In these cases, 
$F(x)~<~x~+~\frac{1}{2}\cdot$
\end{itemize}
Moreover, in cases {\em (iv)} and {\em (v)}, $F(x)$ is the unique real number in $(0,1)$ whose 
base $2$ expansion is given by a periodic balanced sequence and which satisfies 
$x \leq \{F(x) 2^k\} \leq F(x)$ for all $k \geq 0$.
\end{theorem}

\begin{remark}
It is known (see \cite{FM}) that real numbers having a Sturmian base $2$ expansion 
are transcendental. As a consequence of Theorem~\ref{main}, we deduce that, if $x$ is an 
algebraic real number in $[0, \frac{1}{2})$, then $F(x)$ is rational.
\end{remark}

\section{The combinatorial approach}

The main tool used by Bugeaud and Dubickas is {\em combinatorics on words}: real numbers
are replaced by their base $b$ expansion, and inequalities between real numbers
are transformed into (lexicographic) inequalities between infinite sequences representing 
their base $b$ expansions. We will establish a theorem of combinatorial flavour 
(Theorem~\ref{T:main2}), whose translation into a number-theoretical statement is exactly 
Theorem~\ref{main} above. A method for computing $F(x)$ in Case~(v) of Theorem~\ref{main} is 
given in Section~\ref{S:w-from-u}.

\subsection{Two combinatorial theorems}

It happens that the case $b=2$ of Bugeaud-Dubickas' theorem was already proved by 
Veerman in \cite{pV86symb, pV87symb}. The combinatorial result proved by Veerman, 
and by Bugeaud-Dubickas is stated (and strengthened) in Theorems~\ref{P:JPA} and
\ref{T:intro} below. 

\begin{theorem} \label{P:JPA}
An aperiodic sequence $\bs := (s_n)_{n \geq 0}$ on $\{0,1\}$ is Sturmian if
and only if there exists a sequence $\bu := (u_n)_{n \geq 0}$ on $\{0,1\}$
such that $0\bu \leq T^k(\bs) \leq 1\bu$ for all $k \geq 0$. Moreover, $\bu$ is
the unique characteristic Sturmian sequence with the same slope as $\bs$, and we
have $0 \bu = \inf\{T^k(\bs), \ k \geq 0\}$ and $1 \bu = \sup\{T^k(\bs), \ k \geq 0\}$.
\end{theorem}

\begin{theorem} \label{T:intro}
An aperiodic sequence $\bu$ on $\{0,1\}$ is a characteristic
Sturmian sequence if and only if, for all $k \geq 0$,
$$
0\bu < T^k(\bu) < 1\bu.
$$
Furthermore, we have $0 \bu = \inf\{T^k(\bu), \ k \geq 0\}$ and
$1 \bu = \sup\{T^k(\bu), \ k \geq 0\}$.
\end{theorem}

\subsection{The lexicographic world}\label{lex}

As discussed in \cite{AllGlen}, the results in Theorems~\ref{P:JPA} and \ref{T:intro}
have been rediscovered several times since the work of Veerman in the mid-late 80's. 
One of the presentations of these statements is due to Gan \cite{sG01stur}. It is based 
on the {\em lexicographic(al) world}, which seems to have been introduced in 2000, in a 
preprint version of \cite{LabMor2}.

\bigskip

For any two sequences $\bx$, $\by \in \{0,1\}^\NN$, define the set
$$
\varSigma_{\bx,\by} := \{\bs \in \{0,1\}^\NN, \ \forall k \geq 0, \
\bx \leq T^k(\bs) \leq \by\},
$$ 
where $\leq$ denotes the lexicographic order on $\{0,1\}^\NN$ induced by $0 < 1$. 
The {\em lexicographic world} $\mathcal{L}$ is defined by
$$
\mathcal{L} := \{(\mathbf{x},\mathbf{y}) \in \{0,1\}^\NN \times
\{0,1\}^\NN, \, \varSigma_{\mathbf{x,y}} \neq \emptyset\}.
$$
Moreover, by \cite[Lemma~2.1]{sG01stur}, we have 
$$
\cL = \{(\bu,\bv) \in \{0,1\}^\NN \times \{0,1\}^\NN, \, \bv \geq \phi(\bu)\},
$$ 
where $\phi : \{0,1\}^\NN \rightarrow \{0,1\}^\NN$ is the map defined by 
\[
  \phi(\bx) := \inf\{\by \in \{0,1\}^\NN, \, \varSigma_{\bx,\by} \ne \emptyset\}.
\]
Trivially, $\phi(1\bx) = 1^\infty = 111\cdots$ for any sequence $\bx \in \{0,1\}^\NN$.

In \cite{sG01stur}, Gan showed that for any sequence $\bu \in \{0,1\}^\NN$, 
the set $\varSigma_{0\bu, 1\bu}$ is not empty, i.e., there exists a sequence 
$\bs \in \{0,1\}^\NN$ such that $0\bu \leq T^k(\bs) \leq 1\bu$ for all 
$k\geq 0$  (see \cite[Lemma 4.2]{sG01stur}). Furthermore, the sequence $\phi(0\bu)$ has the 
foregoing property (by \cite[Theorem 3.4]{sG01stur}) and it is a Sturmian or 
periodic balanced sequence with the property that $T^k(\phi(0\bu)) \leq \phi(0\bu)$ 
for all $k\geq 0$ (see \cite[Theorem 4.6]{sG01stur}). Moreover,  by \cite[Lemma 5.4]{sG01stur}, 
the set $\varSigma_{0\bu, 1\bu}$ contains a unique Sturmian or periodic balanced sequence 
satisfying $T^k(\bs) \leq \bs$ for all $k \geq 1$. We deduce from these remarks that, for any 
sequence $\bs \in \{0,1\}^\NN$, if $\bs = \phi(0\bu)$ for some sequence $\bu \in \{0,1\}^\NN$, 
then $\bs$ is the unique Sturmian or periodic balanced sequence satisfying 
$0\bu \leq T^k(\bs) \leq 1\bu$ and $T^k(\bs) \leq \bs$ for all $k\geq 0$. The converse of 
this statement also holds by \cite[Corollary 5.6]{sG01stur}. These observations establish Gan's 
main theorem (see below), which shows in particular that any element in the image of 
$\phi$ is a Sturmian or periodic balanced sequence in $\{0,1\}^\NN$ (and such 
sequences are the lexicographically greatest amongst their shifts).

\begin{theorem} \label{T:Gan} \cite[Theorem 1.1]{sG01stur}
For any sequence $\bs \in \{0,1\}^\NN$, the following conditions are equivalent.
\begin{itemize}
\item[(i)] $\bs = \phi(0\bu)$ for some sequence $\bu \in \{0,1\}^\NN$.
\item[(ii)] $\bs$ is the unique Sturmian or periodic balanced sequence satisfying 
$0\bu \leq T^k(\bs) \leq 1\bu$ and $T^k(\bs) \leq \bs$ for all $k \geq 0$.
\end{itemize}
\end{theorem}
\begin{note}
``Sturmian'' in Gan's paper corresponds to what is called here (and classically) 
``Sturmian or periodic balanced''.
\end{note}

\begin{remark}\label{refine} 
It is well known that the closure of the {\em shift-orbit} of a characteristic Sturmian 
sequence $\bs$ (i.e., the closure of $\{T^k(\bs), k\geq 0\}$, denoted by $\Closure(\bs)$) 
is precisely the set of all Sturmian sequences having the same slope as $\bs$ (see for instance 
\cite[Propositions~2.1.25 and 2.1.18]{mL02alge}, or \cite{fM89infi}). In view of this fact, 
Gan's result can be strengthened using Theorem~\ref{P:JPA}, as follows. If the sequence 
$\bs := \phi(0\bu)$ is Sturmian, then $\bu$ is the unique characteristic Sturmian sequence in 
$\Closure(\bs)$, in which case $\bs = 1\bu$. 
\end{remark}

We will further strengthen Gan's result by describing $\phi(0\bu)$ for any given sequence 
$\bu \in \{0,1\}^\NN$. In particular, we will show that when $\bu$ contains both $0$'s and 
$1$'s and is not a characteristic Sturmian sequence, there exists a unique pair of 
characteristic periodic balanced sequences $\bs$ and $\bs'$ of (rational) slope 
$p/q \in (0,1)$ with $\gcd(p,q) = 1$, such that $\bs' \leq \bu \leq \bs$, in which case 
$\phi(0\bu) = 1\bs'$. 
Moreover, the sequences $\bs$, $\bs'$ can be explicitly determined in terms of $\bu$. 

With the same notation as in the Introduction, our main combinatorial theorem reads as follows.

\begin{theorem} \label{T:main2}
Let $\bu$ be a sequence in $\{0,1\}^\NN$. 
\begin{itemize}

\item[(i)] $\phi(1\bu) = 1^\infty$.

\item[(ii)] If $\bu \in \{0^\infty, 1^\infty\}$, then $\phi(0\bu) = \bu$.

\item[(iii)] If $\bu$ is a characteristic Sturmian sequence, then $\phi(0\bu) = 1\bu$. 
Furthermore, $1\bu$ is the unique Sturmian sequence in $\{0,1\}^\NN$ satisfying 
$0\bu \leq T^k(1\bu) \leq 1\bu$ for all $k \geq 0$.

\item[(iv)] If $\bu$ is a characteristic periodic balanced sequence of rational slope 
$p/q \in (0,1)$ with $\gcd(p,q) = 1$, then $\phi(0\bu) = \bs'_{p/q,0} = (1w_{p,q}0)^\infty$.

\item[(v)] If $\bu$ does not take any of the forms given in parts {\em (ii)--(iv)}, then 
there exists a unique pair of coprime integers $p$, $q$ with $0 < p < q$ such that 
$(w_{p,q}01)^\infty < \bu < (w_{p,q}10)^\infty$, in which case $\phi(0\bu) = \bs'_{p/q,0} 
= (1w_{p,q}0)^\infty$.
\end{itemize}
Moreover, in cases (iv) and (v), $\phi(0\bu)$ is the unique periodic balanced sequence in 
$\{0,1\}^\NN$ satisfying $0\bu \leq T^k(\phi(0\bu)) \leq \phi(0\bu)$ for all $k \geq 0$.
\end{theorem}

In the next section, we will recall some generalities about Sturmian and periodic 
balanced sequences. (For more on Sturmian sequences, the reader can consult, e.g., 
\cite[Chapter 2]{mL02alge}.)
The proof of Theorem~\ref{T:main2} is given in Section~\ref{S:main}, and a corollary is 
stated in Section~\ref{S:corollary}. Lastly, in Section~\ref{S:w-from-u}, we show how to 
determine the ``central word'' $w_{p,q}$ such that $\phi(0\bu) = (1w_{p,q}0)^\infty$ for 
any ``generic'' sequence $\bu$ falling into Case~(v) of Theorem~\ref{T:main2} above.

\section{Sturmian \& periodic balanced sequences} \label{S:Sturmian}

In what follows, we will use the following notation and terminology from combinatorics on 
words (see, e.g., \cite{mL02alge}). Let $w = x_1x_2\cdots x_m$ be a word over a finite 
non-empty alphabet $\cA$ (where each $x_i$ is a letter in $\cA$). The \emph{length} of $w$, 
denoted by $|w|$, is equal to $m$. The {\em empty word} is the unique word of length $0$, 
denoted by $\varepsilon$. The number of occurrences of a letter $x$ in $w$ is denoted by 
$|w|_x$. The {\em reversal} of $w$ is defined by $\tilde w = x_m\cdots x_2x_1$, and by 
convention $\varepsilon = \tilde \varepsilon$. If $w = \tilde w$, then $w$ is called a 
{\em palindrome}.  An integer $\ell \geq 1$ is said to be a {\em period} of $w$ if, for 
all $i$, $j$ with $1\leq i, j \leq m$, $i \equiv j \pmod{\ell}$ implies $x_i = x_{j}$. 
Note that any integer $\ell \geq |w|$ is a period of $w$ with this definition. The word 
$w$ is said to be {\em primitive} if it is not a power of a shorter word, i.e., if 
$w = u^n$ implies $n = 1$. 
A finite word $z$ is said to be a {\em factor} of $w$ if $z = x_i x_{i+1} \cdots x_j$ 
for some $i, j$ with $1\leq i \leq j \leq m$. Similarly, a {\em factor} of a sequence 
$\bs := s_0s_1s_2s_3\cdots$ is any finite word of the form $s_is_{i+1}\cdots s_j$ with 
$i \leq j$. 

Recall from the Introduction that the {\em shift map} $T$ is defined on sequences as follows: 
if $\bs := (s_n)_{n\geq 0}$ then $T(\bs) = T((s_n)_{n\geq0}) := (s_{n+1})_{n\geq 0}$. This 
operator naturally extends to finite words as a {\em circular shift} by defining $T(xw) = wx$ 
for any letter $x$ and finite word $w$. 

Under the operation of concatenation, the set $\cA^*$ of all finite words over $\cA$ is a 
{\em free monoid} with identity element $\empt$ and set of generators $\cA$. If $x$ is a 
letter, then we use $x^*$ to denote $\{x\}^*$, the set of all finite powers of $x$. From 
now on, all words and sequences will be over the alphabet $\{0,1\}$.

\subsection{Balanced sequences} \label{SS:balanced}

All Sturmian sequences are ``balanced'' in the following sense (see for instance 
\cite{gHmM40symb, eCgH73sequ, jBpS93acha, jBpS94arem, mL02alge}).

\begin{definition} \label{D:balance} 
A finite word or sequence $w$ over $\{0,1\}$ is said to be {\em balanced} if, 
for any two factors $u$, $v$ of $w$ with $|u| = |v|$, we have $||u|_{1} - |v|_{1}| \leq 1$ 
$($or equivalently $||u|_{0} - |v|_{0}| \leq 1)$. 
\end{definition}

Recall from Remark~\ref{R:slope} that Sturmian sequences are aperiodic. Morse, Hedlund, and 
Coven~\cite{gHmM40symb, eCgH73sequ} proved that the Sturmian sequences are {\em precisely} 
the aperiodic balanced sequences on two letters  (also see \cite[Theorem~2.1.3]{mL02alge}). 
``Periodic balanced sequences'' (as specified in Definition~\ref{D:intro}) are also balanced 
in the sense of the above definition (which justifies their name); moreover, they constitute 
the set of {\em all}  periodic balanced sequences on two letters (see 
\cite[Lemma~2.1.15]{mL02alge} or \cite{rT96onco}). 

\subsection{Characteristic Sturmian sequences} 

In \cite{aD97stur}, characteristic Sturmian sequences were characterized using iterated 
palindromic closure, defined as follows. The {\em palindromic (right-)closure} of a finite 
word $w$, denote by $w^{(+)}$, is the (unique) shortest palindrome beginning with $w$. That is, 
if $w = uv$ where $v$ is the longest palindromic suffix of $w$, then $w^{(+)} := uv\tilde u$.
For example, $(011)^{(+)} = 0110$. The {\em iterated palindromic closure function}, denote by 
$Pal$, is defined by iteration of the  palindromic right-closure operator (see, e.g., 
\cite{jJ05epis}). More precisely, $Pal$ is defined recursively as follows. 
Set $Pal(\empt) = \empt$,  and for any word $w$ and letter $x$, define 
$Pal(wx) := (Pal(w)x)^{(+)}$. For example, $Pal(011) = (Pal(01)1)^{(+)} = (0101)^{(+)} = 01010$. Note that $Pal$ is injective; and moreover, it is clear from the definition that $Pal(w)$ is 
a prefix of $Pal(wx)$ for any word $w$ and letter $x$. Hence, if $v$ is a prefix of $w$, then 
$Pal(v)$ is a prefix of $Pal(w)$.

The following theorem provides a combinatorial description of characteristic Sturmian 
sequences in terms of $Pal$.

\begin{theorem} \label{T:Sturmian} \cite{aD97stur}
For any sequence $\bs \in \{0,1\}^\NN$, the following properties are equivalent.
\begin{itemize}
\item[(i)] $\bs$ is a characteristic Sturmian sequence.
\item[(ii)] There exists a (unique) sequence $\Delta := x_0x_1x_2x_3\ldots \in 
\{0,1\}^\NN \setminus \left(\{0,1\}^*0^\infty \cup \{0,1\}^*1^\infty\right)$ (i.e., 
not eventually constant), called the {\em directive sequence} of $\bs$, such that 
\[ 
\bs = \lim_{n\rightarrow\infty}Pal(x_0x_1x_2\cdots x_n) = Pal(\Delta).
\]
\end{itemize}
\end{theorem}

\begin{example} \label{Ex:Sturmian} Recall from Example~\ref{Ex:Fibonacci} that the (binary) 
Fibonacci sequence $\bbf = 01001010010\cdots$ is the characteristic Sturmian sequence 
$\bc_\alpha$ with $\alpha = (3 - \sqrt{5})/2$; it has directive sequence $(01)^\infty$. That is:
\[
\bbf = Pal(0101\cdots) = \underline{0}\underline{1}0\underline{0}10\underline{1}0010\cdots,
\]
where the underlined letters indicate at which points palindromic closure is applied. Note the 
simple continued fraction expansion of $\alpha = (3 - \sqrt{5})/2$ is $[0;2,1,1,1,\ldots]$. 
More generally, if $\alpha \in (0,1)$ is an irrational number with simple continued fraction 
expansion $[0; d_1+1,d_2,d_3,d_4,\ldots]$ where $d_1 \geq 0$ and all other $d_i \geq 1$, then 
$\bc_\alpha = Pal(0^{d_1}1^{d_2}0^{d_3}1^{d_4}\cdots)$ (see \cite{aFmMuT78dete, tB93desc} and 
also \cite[pg. 206]{pAgR91repr}). 
\end{example}

\subsection{Characteristic periodic balanced sequences} \label{SS:periodic}

We will now recall some known combinatorial descriptions of the characteristic 
periodic balanced sequences in $\{0,1\}^\NN$ (see Proposition~\ref{P:central}  
and Remark~\ref{R:char-periodic} below).

Let us first recall from the Introduction that the characteristic balanced sequences of slopes 
$0$ and $1$ are $\bc_0 = 0^\infty$ and $\bc_1 = 1^\infty$, respectively. For all other rational 
slopes $p/q \in (0,1)$ with $\gcd(p,q) = 1$, there are exactly two characteristic periodic 
balanced sequences of slope $p/q$, given by
\[
\bs_{p/q,p/q} = T(\bs_{p/q,0}) = (w_{p,q}10)^\infty \quad \mbox{and} \quad \bs'_{p/q,p/q} 
= T(\bs'_{p/q,0}) = (w_{p,q}01)^\infty,
\]
where $w_{p,q}$ is a word of length $q - 2$ in $\{0,1\}^*$. For example, with $p = 2$ and 
$q = 5$, we obtain the following two characteristic periodic balanced sequences of slope $2/5$:
\[
 \bs_{2/5,2/5} = (01010)^\infty \quad \mbox{and} \quad \bs'_{2/5,2/5} = (01001)^\infty \quad 
\mbox{where $w_{2,5} = 010$.}
 \]
Notice that $w_{2,5}$ is a palindrome and $|w_{2,5}10|_1 = |w_{2,5}01|_1 = 2 = p$. 
More generally, one can verify that all words $w_{p,q}$ are palindromes and 
$|w_{p,q}10|_{1} = |w_{p,q}01|_1 = p$. Furthermore, the words $w_{p,q}10$ and $w_{p,q}01$ 
(which have length $q$) are primitive since $\gcd(p,q) = 1$. Hereafter, the word $w_{p,q}$ 
will be called the {\em central word of slope $p/q$}; it is the unique central word of length 
$q - 2$ containing $p - 1$ occurrences of $1$.

\begin{note} The set of all central words of slope $p/q \in (0,1)$ (where $p$, $q$ are 
coprime integers) coincides with the family of ``central words'' in $\{0,1\}^*$ as defined 
in \cite[Chapter 2]{mL02alge} (in particular, see 
\cite[Theorem~2.2.11 and Proposition~2.2.12]{mL02alge}).
\end{note}

The following proposition collects together some equivalent definitions of central words. 
For many more, see the nice survey~\cite{jB07stur}.

\begin{proposition} \label{P:central}
For any word $w \in \{0,1\}^*$, the following properties are equivalent.
\begin{itemize}
\item[(i)] $w$ is a central word.
\item[(ii)] $0w1$ and $1w0$ are balanced \cite{aDfM94some}. 
\item[(iii)] $w = Pal(v)$ for some word $v \in \{0,1\}^*$ \cite{aDfM94some, aD97stur}.
\item[(iv)] $w$ has two periods $\ell$, $m$ such that $\gcd(\ell,m) = 1$ and $|w|=\ell+m-2$ 
\cite{eCgH73sequ, aDfM94some}.
\item[(v)] $w \in 0^* \cup 1^* \cup (P \cap P10P)$ where $P$ is the set of all palindromes 
in $\{0,1\}^*$ \cite{aDfM94some}.
\item[(vi)] $w \in 0^*\cup 1^*$, or there exists a unique pair of words 
$w_1$, $w_2 \in \{0,1\}^*$ such that $w$ satisfies the equation $w = w_101w_2 = w_2 10 w_1$ 
\cite{aDfM94some, aD97stur}. 
\end{itemize}
Moreover, in part {\em (vi)}, $w_1$ and $w_2$ are central words, $\ell_1 := |w_1| + 2$ and 
$\ell_2 := |w_2| + 2$ are coprime periods of $w$, and $\min\{\ell_1, \ell_2\}$ is the minimal 
period of $w$ \cite{aCaD05code}.
\end{proposition}
\begin{note} $P \cap (P10P) = P \cap (P01P)$.
\end{note}

Furthermore, by \cite[Proposition~2.2.12]{mL02alge}, the central word $w_{p,q}$ of slope 
$p/q$ is the central word with coprime periods $\ell$, $m$ where $\ell + m = q$ and 
$m p \equiv 1 \pmod{q}$. For example, the central word $w_{2,5} = 010$ has coprime periods 
$\ell = 2$ and $m = 3$ where $2 + 3 = q$ and $m p = 6 \equiv 1 \pmod{5}$. Also note that 
$w_{1,q} = 0^{q-2}$ and $w_{q-1,q} = 1^{q-2}$; in particular $w_{1,2} = \varepsilon$.

\begin{remark} \label{R:char-periodic}
Let $p$, $q$ be coprime integers with $0 < p < q$. Then $p/q$ has two distinct simple 
continued fraction expansions:
\[
p/q = [0; d_1 + 1, \ldots, d_n, 1] = [0; d_1 + 1, \ldots, d_n + 1]
\]
where $d_1 \geq 0$ and all other $d_i \geq 1$. It is known (see, e.g., 
\cite[Proposition~27]{jB07stur}) that the word $v \in \{0,1\}^*$ such that 
$w_{p,q} = Pal(v)$ takes the form $v = 0^{d_1}1^{d_2}0^{d_3}\cdots x^{d_n}$ 
where $x = 0$ if $n$ is odd and $x = 1$ if $n$ is even. For example, 
$2/5 = [0; 2,1,1] = [0; 2,2]$ and $w_{2,5} = 010 = Pal(01)$. Moreover, as in 
the case of characteristic Sturmian sequences (see Theorem~\ref{T:Sturmian} 
and Example~\ref{Ex:Sturmian}), the two characteristic periodic balanced sequences 
of slope $p/q$ can be obtained by iterated palindromic closure. More precisely, with 
the above notation, we have
\[
(w_{p,q}xy)^\infty = Pal(0^{d_1}1^{d_2} \cdots x^{d_n+1}y^\infty) \quad \mbox{and} \quad 
(w_{p,q}yx)^\infty = Pal(0^{d_1}1^{d_2}\cdots x^{d_n}yx^\infty)
\]
where $\{x,y\} = \{0,1\}$. For example, 
\[
(w_{2,5}10)^\infty = (01010)^\infty = Pal(0110^\infty) \quad \mbox{and} 
\quad (w_{2,5}01)^\infty = (01001)^\infty = Pal(0101^\infty).
\]
\end{remark}

\medskip
In \cite{gP99anew} Pirillo proved that a word $w \in \{0,1\}^*$ is a palindromic prefix 
of some characteristic Sturmian sequence in $\{0,1\}^\NN$, i.e., $w = Pal(v)$ for some 
$v \in \{0,1\}^*$ (see  Theorem~\ref{T:Sturmian}) if and only if $w01$ is a circular shift 
of $w10$. From this fact and Proposition~\ref{P:central}, we thus deduce the following result. 

\begin{proposition} \label{P:circular-shift}
A word $w \in \{0,1\}^*$ is central if and only if $w01$ is a circular shift of $w10$.
\end{proposition}

Consequently, for any two coprime integers $p$, $q$ with $0 < p < q$, the two characteristic 
periodic balanced sequences of slope $p/q$, namely $\bs_{p/q,p/q} = (w_{p,q}10)^\infty$ and 
$\bs'_{p/q,p/q} = (w_{p,q}01)^\infty$, are shifts of each other, and therefore they have the 
same set of factors.

\begin{remark}\label{R:orbit} Recall from Remark~\ref{refine} that the closure of the 
shift-orbit of a characteristic Sturmian sequence $\bs$ is precisely the set of all 
Sturmian sequences having the same slope as $\bs$. Moreover, the (Sturmian) sequences in 
$\Closure(\bs)$ are exactly the sequences that have the same set of factors as $\bs$ 
(see for instance \cite[Propositions~2.1.25 and 2.1.18]{mL02alge}, or \cite{fM89infi}). 
Likewise, the shift-orbit of a characteristic periodic balanced sequence $\bu$ consists of 
all the periodic balanced sequences with the same set of factors (and also the same slope) 
as $\bu$. However, in contrast to the aperiodic case, we deduce from 
Proposition~\ref{P:circular-shift} that, if $\bu$ is a characteristic periodic balanced 
sequence containing both $0$'s and $1$'s, then $\Closure(\bu)$ contains exactly two distinct 
characteristic periodic balanced sequences (not just one), which take the form $(w01)^\infty$ 
and $(w10)^\infty$ where $w$ is the central word having the same slope as $\bu$.
\end{remark}

The following useful result is due to de~Luca~\cite{aD97stur}; in particular, 
see \cite[Remark~1 and Proposition~9]{aD97stur} and also \cite[Lemma~5]{aCaD05code}.

\begin{proposition} \label{P:deLuca} Let $w$ be a central word in $\{0,1\}^*$. 
If $w = w_101w_2 = w_210w_1$ where $w_1$ and $w_2$ are (central) words, then 
\[
Pal(w0) = (w0)^{(+)} = w_210w_101w_2 \quad \mbox{and} 
\quad Pal(w1) = (w1)^{(+)} = w_101w_210w_1.
\]
\end{proposition}

We end this section with a result (Corollary~\ref{C:Berstel-deLuca} below) that will be 
particularly useful in the proof of our main combinatorial theorem. Let us first recall 
that a non-empty finite word $v$ over an alphabet $\cA$ is said to be a {\em Lyndon word} 
(resp.~{\em anti-Lyndon word}) if $v$ is a primitive word that is  lexicographically less 
(resp.~lexicographically greater) than all of its circular shifts with respect to a given 
total order on~$\cA$.

\begin{proposition} \label{P:Berstel-deLuca} \cite[Theorem~3.2 and Corollary~3.1]{jBaD97stur}
A non-empty finite word $v \in \{0,1\}^*$ is a balanced Lyndon word (resp.~balanced anti-Lyndon 
word) with respect to the lexicographic order if and only if  $v = 0w1$ (resp.~$v = 1w0$) for 
some central word $w \in \{0,1\}^*$.
\end{proposition}

As a direct consequence of the above proposition, we have the following result.

\begin{corollary} \label{C:Berstel-deLuca}
For any central word $w \in \{0,1\}^*$, the (primitive) words $0w1$ and $1w0$ are the 
lexicographically least and greatest words amongst their circular shifts.
\end{corollary}

\section{Proof of Theorem \ref{T:main2}} \label{S:main}

-- Assertions (i) and (ii) are straightforward (see \cite[Lemma~2.4]{sG01stur}).

\medskip

In order to prove the other assertions, we first recall the inequalities that are equivalent 
to $\bs = \phi(0\bu)$ from Theorem~\ref{T:Gan}: 
\begin{equation} \label{eq:inequalities}
0\bu \leq T^k(\bs) \leq 1\bu \quad \mbox{and} \quad T^k(\bs) \leq \bs \quad 
\mbox{for all $k \geq 0$}.
\end{equation}

-- Assertion (iii) is a consequence of Gan's result (Theorem~\ref{T:Gan}) together with 
Theorem~\ref{P:JPA} (see Remark~\ref{refine}).
  
\medskip

-- We will now prove Assertions (iv) and (v). Suppose that $\bu$ is not a characteristic 
Sturmian sequence and that $\bu$ contains both $0$'s and $1$'s. Then we know from 
Theorem~\ref{T:Gan} that $\bs := \phi(0\bu)$ is a periodic balanced sequence satisfying the 
inequalities in \eqref{eq:inequalities}. Indeed, $\bs$ cannot be Sturmian, for otherwise $\bu$ 
would be a characteristic Sturmian sequence by Remark~\ref{refine}. 
   
By Remark~\ref{R:orbit}, $\Closure(\bs)$ contains exactly two distinct characteristic periodic 
balanced sequences, given by
\[
\bs_{01} :=  (w01)^\infty \quad \mbox{and} \quad \bs_{10} := (w10)^\infty
\]
where $w \in \{0,1\}^*$ is the central word with the same slope as $\bs$.  
We now deduce from Corollary~\ref{C:Berstel-deLuca} that the lexicographically least 
sequence in $\Closure(\bs)$ is
\[
(0w1)^\infty = 0(w10)^\infty = 0\bs_{10} 
\]
and the lexicographically greatest sequence in $\Closure(\bs)$ is
\[
(1w0)^\infty = 1(w01)^\infty = 1\bs_{01}.
\]
Hence, 
\begin{equation} \label{eq;01}
0\bs_{10} \leq T^k(\bs) \leq 1\bs_{01} \quad \mbox{for all $k \geq 0$}.
\end{equation}
Moreover, since $\bs$ is the lexicographically greatest sequence in its shift-orbit 
(by the second  inequality in \eqref{eq:inequalities}), we have 
$\bs = (1w0)^\infty = 1\bs_{01}$. 

We will now show that $\bs_{01} \leq \bu \leq \bs_{10}$. Since $0\bs_{10}$ and $1\bs_{01}$ 
are the lexicographically least and greatest elements in $\Closure(\bs)$, the inequalities 
in \eqref{eq:inequalities} imply that 
\[
0 \bu \leq 0\bs_{10} \quad \mbox{\rm and} \quad 1\bs_{01} \leq 1\bu.
\]
Hence $\bs_{01} \leq \bu \leq \bs_{10}$; that is, $(w01)^\infty \leq \bu \leq (w10)^\infty$.

Furthermore, we note that there does not exist another central word $z$ 
such that $(z01)^\infty \leq \bu \leq (z10)^\infty$. For if so, then the set 
$[0\bu,1\bu] := \{\bs \in \{0,1\}^\NN, 0\bu \leq \bs \leq 1\bu\}$ would contain the periodic 
balanced sequences $0(z10)^\infty = (0z1)^\infty$ and $1(z01)^\infty = (1z0)^\infty$, and 
hence all of the shifts of the characteristic periodic balanced sequence $(z01)^\infty$, 
since the former two sequences are the lexicographically least and greatest sequences in 
the shift-orbit of $(z01)^\infty$ (by Proposition~\ref{P:circular-shift} and 
Corollary~\ref{C:Berstel-deLuca}). But by \cite[Lemma 5.4]{sG01stur}, the set $[0\bu,1\bu]$ 
contains a unique periodic balanced shift-orbit. Therefore, since 
$\Closure((w01)^\infty) \subseteq [0\bu,1\bu]$, we must have $z = w$. We have thus established 
the following lemma.

\begin{lemma} \label{L:u-bounded}
Suppose $\bu$ is a sequence in $\{0,1\}^\NN \setminus\{0^\infty, 1^\infty\}$ that is not 
characteristic Sturmian. Then there exists a unique central word $w \in \{0,1\}^*$ such that 
$(w01)^\infty \leq \bu \leq (w10)^\infty$. Moreover, $\phi(0\bu) = (1w0)^\infty$.
\end{lemma}

Assertions (iv) and (v) are direct consequences of the above lemma, and the last statement 
in the theorem follows from Theorem~\ref{T:Gan}. \qed

\section{A corollary} \label{S:corollary}

By Theorem~\ref{T:Sturmian}, the set of all characteristic Sturmian sequences 
in $\{0,1\}^\NN$ is given by
\[
\cS = \{\bu \in \{0,1\}^\NN, \, \exists \bv \in \{0,1\}^\NN \setminus \left(\{0,1\}^*0^\infty 
\cup \{0,1\}^*1^\infty\right), \, \bu = Pal(\bv)\}.
\]
And it follows from Proposition~\ref{P:central} that the set of all characteristic periodic 
balanced sequences in $\{0,1\}^\NN$ is given by
\[
\cP = \{0^\infty, 1^\infty\} \cup \cP_{01} \cup \cP_{10} 
\]
where
\[
\cP_{01} := \{\bu \in \{0,1\}^\NN, \, \exists v \in \{0,1\}^*, \, 
\bu = (Pal(v)01)^\infty\}
\]
and
\[
\cP_{10} :=  \{\bu \in \{0,1\}^\NN, \, \exists v \in \{0,1\}^*, \, 
\bu = (Pal(v)10)^\infty\}.
\]
Given a characteristic periodic balanced sequence in $\{0,1\}^\NN$ of the form 
$\bs := (Pal(v)xy)^\infty$ where $v \in \{0,1\}^*$ and $\{x,y\} = \{0,1\}$, we let $\bar \bs$ 
denote the other characteristic periodic balanced sequence in the shift-orbit of $\bs$, i.e., 
$\bar \bs := (Pal(v)yx)^\infty$ (see Remark~\ref{R:orbit}).

As an immediate consequence of Theorem~\ref{T:main2},   
we obtain the following description of the lexicographic world.

\begin{corollary} \label{C:main} We have 
$\cL = \{(01^\infty, 1^\infty)\} \cup \cL_0 \cup  \cL_1 \cup \cL_{01} \cup \cL_{10} 
\cup \cL_{*}$ 
where 
\begin{eqnarray*}
\cL_0 &:=& \{(0^\infty, \bv), \, \bv \in \{0,1\}^\NN\}, \\
\cL_1 &:=& \{(1\bu,1^\infty), \, \bu \in \{0,1\}^\NN\}, \\ 
\cL_{01} &:=& \{(0{\bu}, \bv) \in 0(\cS\cup\cP_{01}) \times \{0,1\}^\NN, \, 
\bv \geq 1\bu\}, \\ 
\cL_{10} &:=& \{(0\bu,\bv) \in 0\cP_{10} \times \{0,1\}^\NN, \, \bv \geq 1\bar \bu\}, \\
\cL_{*} &:=& \{(0\bu,\bv) \in (\{0,1\}^\NN\setminus 0(\cS \cup \cP)) \times \{0,1\}^\NN, \,  
\exists \bs \in \cP_{01}, \bs \leq \bu \leq \bar\bs, \bv \geq 1\bs\}.
\end{eqnarray*}
\end{corollary}

\section{How to determine $\phi(0\bu)$ for a generic sequence $\bu$} \label{S:w-from-u}

The following theorem provides a method for determining the central word $w$ such that 
$\phi(0\bu) = (1w0)^\infty$ for any ``generic'' sequence $\bu \in \{0,1\}^\NN$ falling into 
Case~(v) of Theorem~\ref{T:main2}. Hereafter, a prefix of $\bu$ that is a central word is 
called a {\em central prefix} of $\bu$. 

\begin{theorem} \label{T:w-from-u}
Suppose $\bu$ is a sequence in $\{0,1\}^\NN\setminus \{0^\infty, 1^\infty\}$ that is neither 
a characteristic Sturmian sequence nor a characteristic periodic balanced sequence. Let $v$ 
be the longest central prefix of $\bu$. Then $v$ is finite ($v\ne \varepsilon$), and 
$\phi(0\bu)$ is determined as follows.
\begin{itemize}
\item[(i)] If $v = 1^k$ for some $k \geq 1$, then 
$\phi(0\bu) =  (1^k0)^\infty = (1w_{p,q}0)^\infty$ where $p = k$ and $q = k + 1$.
\item[(ii)] If $v = 0^k$ for some $k \geq 1$, then 
$\phi(0\bu) =  (10^k)^\infty = (1w_{p,q}0)^\infty$ where $p = 1$ and $q = k + 1$.
\item[(iii)] Suppose $v$ contains both $0$'s and $1$'s. Let $v_1$, $v_2$ be the unique pair 
of central words such that $v = v_101v_2 = v_2 10 v_1$ where $\ell_1 := |v_1| + 2$ and 
$\ell_2 := |v_2| + 2$ are coprime periods of $v$. Consider the prefix of length $2|v| + 4$ 
of $\bu$, namely the prefix $vxyz$ where $x, y \in \{0,1\}$ and $|z| = |v| + 2$.
\begin{itemize}
\item[(a)] If either $xy = 01$ and $z > v01$, or $xy = 10$ and $z < v10$, then 
$\phi(0\bu) = (1v0)^\infty = (1w_{p,q}0)^\infty$ where $p=|v|_1+1$ and $q=|v|+2=\ell_1+\ell_2$. 
\item[(b)] If either $xy = 01$ and  $z < v01$, or $xy = 00$, then 
$\phi(0\bu) = (1v_20)^\infty = (1w_{p,q}0)^\infty$ where $p = |v_2|_1 + 1$ and $q = \ell_2$.
\item[(c)] If either $xy = 10$ and $z > v10$, or $xy = 11$, then 
$\phi(0\bu) = (1v_10)^\infty = (1w_{p,q}0)^\infty$ where $p = |v_1|_1 + 1$ and $q = \ell_1$.
\end{itemize}
\end{itemize}
\end{theorem}
\begin{note} In Assertion~(iii), it cannot happen that $z = vxy$ when $x \ne y$. For instance, 
if $xy = 01$ and $z = v01$, then $\bu$ would begin with the following word:
\[
v01v01 = v_210v_101v_210v_101
\]
where the prefix $v_210v_101v_2 = v01v_2$ is a central word, by Propositions~\ref{P:central} 
and \ref{P:deLuca}. But then $\bu$ has a central prefix longer than $v$; thus $z \ne v01$. 
Similarly, if $xy = 10$, then $z \ne v10$.
\end{note}

The following lemma is needed for the proof of Theorem~\ref{T:w-from-u}.

\begin{lemma} \label{L:w-from-u} Suppose $v$ is a central word in 
$\{0,1\}^*\setminus (0^*\cup 1^*)$. Let $v_1$, $v_2$ be the unique pair of central 
words such that $v$ satisfies the equation $v = v_101v_2 = v_2 10 v_1$. Then
$v01v_2$ (resp.~$v10v_1$) is a prefix of the characteristic periodic balanced 
sequence $(v_210)^\infty$ (resp. $(v_101)^\infty$).
\end{lemma} 
\begin{note} By Propositions~\ref{P:central} and \ref{P:deLuca}, the words $v01v_2$ 
and $v10v_1$ are central words since $v01v_2 = Pal(v0)$ and $v10v_1 = Pal(v1)$.
\end{note}
\begin{proof}[Proof of Lemma~$\ref{L:w-from-u}$]
We will prove only that $(v_210)^\infty$ begins with the central word $Pal(v0) = v01v_2$, 
since the proof of the other case is very similar.

By Proposition~\ref{P:central}, $\ell_1 := |v_1| + 2$ and $\ell_2 := |v_2| + 2$ are coprime 
periods of $v$ where $|v| = \ell_1 + \ell_2 - 2$. In particular, since $\ell_2 = |v_210|$ is 
a period of $v$ with $\ell_2 < |v|$, there exists an integer $k\geq 1$ such that 
$v = (v_210)^kv_2'$ where $v_2'$ is a (possibly empty) prefix of $v_210$, in which 
case $v_1 = (v_210)^{k-1}v_2'$ since $v = v_210v_1$. 
Moreover, since $v$ is a palindrome, $\tilde v_2'$ is a prefix of $v$, 
and therefore $v_2' = \tilde v_2'$, i.e., $v_2'$ is a palindrome. Furthermore, $v_2'01$ is a 
prefix of $v$ since its reversal $10v_2'$ is a suffix of~$v$. We will now show that the central 
word $Pal(v0) = v01v_2$ is a prefix of the characteristic periodic balanced sequence 
$(v_210)^\infty$ by considering five different cases according to the length of the palindrome 
$v_2'$.

\medskip

\begin{enumerate}
\item If $v_2' = v_210$, then $v = (v_210)^{k+1}$ and $v_1 = (v_210)^k$. Hence, $v$ being a 
palindrome implies that $v_210$ is a palindrome and we have $v = (v_210)^{k+1} = (01v_2)^{k+1}$. 
Therefore,
\begin{equation*}
(v_210)^\infty = (v_210)^{k+1}\underbrace{01v_2}_{v_210}(v_210)^\infty = v01v_2(v_210)^\infty.
\end{equation*}
Thus, the central word $Pal(v0) = v01v_2$ is a prefix of $(v_210)^\infty$.

\medskip
\item If $v_2' = v_21$, then since $v_2'$ and $v_2$ are palindromes, we have $v_21 = 1v_2$, 
and hence $v_2$ is a power of $1$; in particular, $v_2 = 1^{\ell_2-2}$. Therefore
\begin{eqnarray*}
(v_210)^\infty &=& (v_210)^kv_210v_210(v_210)^\infty \\
&=& \underbrace{(v_210)^kv_21}_{v}0\underbrace{11^{\ell_2-2}}_{1v_2}0(v_210)^\infty \\
&=& v01v_20(v_210)^\infty \\
&=& Pal(v0)0(v_210)^\infty.
\end{eqnarray*}

\medskip

\item If $v_2' = v_2$, then $v = (v_210)^{k}v_2$, and therefore $v_1 = (v_210)^{k-1}v_2$. 
But this implies that $\ell_1 = k\ell_2$, which is impossible since $\ell_1$ and $\ell_2$ 
are coprime integers greater than $1$.

\medskip
\item If $v_2 = v_2'0$, then since $v_2$ and $v_2'$ are palindromes, we have
$v_2'0 = 0v_2'$. Therefore $v_2'$ (and hence $v_2$) is a power of $0$;
in particular, $v_2 = 0^{\ell_2-2}$. Thus
\begin{eqnarray*}
(v_210)^\infty &=& (v_210)^kv_210v_210(v_210)^\infty \\
&=& \underbrace{(v_210)^kv_2'}_{v}01\underbrace{v_20}_{0v_2}10(v_210)^\infty \\
&=& v01v_2010(v_210)^\infty \\
&=&Pal(v0)010(v_210)^\infty.
\end{eqnarray*}

Note that we cannot have $v_2 = v_2'1$ because $v_2'01$ and $v_2$ are both prefixes of $v$.

\medskip
\item If $|v_2'| \leq |v_2| - 2$, then $v_2 = v_2'01v_2''$ for some (possibly empty) word 
$v_2'' \in \{0,1\}^*$, in which case $v = (v_2'01v_2''10)^kv_2'$.  (Note that neither  
$v_2'1$ nor $v_2'00$ is a prefix of $v_2$ because $v_2'01$ and $v_2$ are both prefixes 
of $v$.) Since $v$ is a palindrome that begins with the palindrome $v_2 = v_2'01v_2''$ 
and therefore ends with $\tilde v_2 = v_2 = \tilde v_2''10v_2'$, we see that 
$\tilde v_2''10v_2' = v_2''10v_2'$. Hence $v_2''$ is a palindrome. Moreover, $v_2''$ is 
a central word since $v_2''$ is a palindromic prefix (and also a palindromic suffix) of 
the central word $v_2$ and any palindromic prefix (or suffix) of a central word is central 
(see \cite{aDfM94some} or \cite[Corollary 2.2.10]{mL02alge}). Thus, by 
Proposition~\ref{P:central}, $v_2$ satisfies the equation $v_2 = v_2''10v_2' = v_2'01v_2''$.  
Hence, we have
\begin{eqnarray*}
(v_210)^\infty 
&=& (v_210)^k\underbrace{v_2'01v_2''10}_{v_210}v_210(v_210)^\infty  \\ 
&=& (v_210)^kv_2'01\underbrace{v_2''10v_2'01}_{v_201}v_2''10(v_210)^\infty   \\
&=& v01v_201v_2''10(v_210)^\infty \\
&=& Pal(v0)01v_2''10(v_210)^\infty.
\end{eqnarray*}
\end{enumerate}

In all of the above cases (with the exception of the impossible case~(3)), we have shown 
that the central word $Pal(v0) = v01v_2$ is a prefix of $(v_210)^\infty$, as required.
\end{proof}

\begin{proof}[Proof of Theorem~$\ref{T:w-from-u}$] Suppose $\bu$ is a sequence in 
$\{0,1\}^\NN \setminus \{0^\infty, 1^\infty\}$ that is neither a characteristic Sturmian 
sequence nor a characteristic periodic balanced sequence. Then the longest central prefix 
of $\bu$, say $v$, is non-empty since it could (at the very least) be a letter. Furthermore, 
$v$ is finite; otherwise, if $v$ were infinite, then $\bu$ would be either a characteristic 
Sturmian sequence or a characteristic periodic balanced sequence (see Theorem~\ref{T:Sturmian} 
and Remark~\ref{R:char-periodic}). 

We know from Theorem~\ref{T:main2}  (or Lemma~\ref{L:u-bounded}) that there exists a unique 
central word $w \in \{0,1\}^*$ such that $(w01)^\infty < \bu < (w10)^\infty$, in which case 
$\phi(0\bu)$ is equal to the periodic balanced sequence $(1w0)^\infty$. We will show how to 
determine $w$ in terms of the longest central prefix $v$. Note that $w$ is either empty or a 
(palindromic) prefix of $v$, by the maximality of $v$. 

First suppose that $v = x^k$ for some $x \in \{0,1\}$ and $k \geq 1$. Then by the maximality 
of $v$ as a central prefix of $\bu$, it follows that $\bu$ begins with $x^ky = vy$ where 
$y \in \{0,1\}$, $y \ne x$. Moreover, the prefix of length $2k +1$ of $\bu$ takes the form 
$x^kyu$ where $|u| = k$ and $|u|_x \leq k-1$; otherwise $\bu$ would begin with 
$x^kyx^k = Pal(x^ky)$, contradicting the fact that $v ~(= Pal(x^k))$ is the longest central 
prefix of $\bu$. If $x = 1$, then we easily see that 
\[
(1^{k-1}01)^\infty < \bu~(= 1^k0u \cdots) < (1^k0)^\infty,
\]
where the latter inequality follows from the fact that $|u| = k$ and $u < 1^k$ (since $u$ 
contains at most $k-1$ occurrences of the letter $1$). Hence by Lemma~\ref{L:u-bounded}, 
$\phi(0\bu) =  (1^k0)^\infty = (1w_{p,q}0)^\infty$ where $p = k$ and $q = k + 1$. Similarly, 
if $x = 0$, we have
\[
(0^k1)^\infty < \bu~(= 0^k1u \cdots) < (0^{k-1}10)^\infty,
\]
where the first inequality follows from the fact that $|u| = k$ and $0^k < u$ (since $u$ 
contains at most $k-1$ occurrences of the letter $0$). Hence by Lemma~\ref{L:u-bounded}, 
$\phi(0\bu) =  (10^k)^\infty = (1w_{p,q}0)^\infty$ where $p = 1$ and $q = k + 1$. We have 
thus proved Assertions (i) and (ii) of the theorem.

Now suppose that the longest central prefix $v$ of $\bu$ contains both 0's and 1's. 
Then by Proposition~\ref{P:central}, there exists a unique pair of central words 
$v_1$, $v_2 \in \{0,1\}^*$ such that $v = v_101v_2 = v_2 10 v_1$ where $\ell_1 := |v_1| + 2$ 
and $\ell_2 := |v_2| + 2$ are coprime periods of $v$,  and $\min\{\ell_1, \ell_2\}$ is the 
minimal period of~$v$. 

Consider the prefix of length $2|v| + 4$ of $\bu$, namely the prefix $vxyz$ where 
$x, y \in \{0,1\}$ and $|z| = |v| + 2$. We will now prove each of the cases (a), (b), and (c) 
of Assertion~(iii).

\medskip

\begin{description}
\item[Case (a)] Let us first suppose that $\bu$ begins with $v01z$ where $|z| = |v01|$ and 
$z > v01$. Then it is easy to see that
\[
 (v01)^\infty < \bu < (v10)^\infty.
\]
Hence, by Lemma~\ref{L:u-bounded}, we have $\phi(0\bu) = (1v0)^\infty$. Moreover, 
$v = w_{p,q}$ where $p = |v|_1 + 1$ and $q = |v| + 2 = \ell_1 + \ell_2$.  Similarly, 
if $\bu$ begins with $v10z$ where $|z| = |v10|$ and $z < v10$, then 
$(v01)^\infty < \bu < (v10)^\infty$, and therefore $\phi(0\bu) = (1v0)^\infty$ by 
Lemma~\ref{L:u-bounded}. 

\medskip

\item[Case (b)] In this case, either $\bu$ begins with $v00$ or $\bu$ begins with $v01z$ 
where $|z| = |v01|$ and $z < v01$. 

Since $v_210$ is a prefix of $v$ (which in turn is a prefix of $\bu$), we have 
$(v_201)^\infty < \bu$. Furthermore, by Lemma~\ref{L:w-from-u}, the characteristic 
periodic balanced sequence $(v_210)^\infty$ begins with the central word $Pal(v0) = v01v_2$. 
Thus, if $\bu$ begins with $v00$, then $\bu < (v_210)^\infty$. 
On the other hand, if $\bu$ begins with $v01z$ where $|z| = |v01|$ and $z < v01$, then we 
will show that $\bu < (v_210)^\infty$ by considering the prefix of length $|v| + \ell_2$ of 
$\bu$, namely $v01z_2$ where $|z_2| = |v_2|$. We first note that $z_2 \leq v_2$ since $v_2$ 
is a prefix of $v$ and $z_2$ is a prefix of $z$ where $z$ and $v$ satisfy $z < v01$. 
Furthermore, $z_2 < v_2$ (i.e., $z_2 \ne v_2$). Otherwise, if $z_2 = v_2$, then $\bu$ would 
begin with the central word $Pal(v0) = v01v_2$. But then $\bu$ would have a central prefix 
that is longer than $v$; a contradiction. Therefore $z_2 < v_2$, and hence 
$\bu < (v_210)^\infty$ since $(v_210)^\infty$ begins with $v01v_2$ where $v_2 > z_2$, as 
shown above.

\medskip

\item[Case (c)] This case is symmetric to Case (b).
\end{description}
\end{proof}

\begin{example} \label{Ex:computation}
The following examples demonstrate the computation of $\phi(0\bu)$ for sequences 
$\bu$ in $\{0,1\}^\NN$ that are neither characteristic Sturmian nor periodic and balanced. 
Where appropriate, the longest central prefix of the sequence is highlighted in boldface. 
\begin{enumerate}
\item The following two general facts can easily be deduced from the proofs of parts (i) 
and (ii) of Theorem~\ref{T:w-from-u}.

\smallskip

\begin{enumerate}
\item $\phi(0\bu) = (1^k0)^\infty$ for any sequence $\bu$ having a prefix of the form $1^k0v$ 
where $k\geq 1$, $|v| = k$, and $|v|_1 \leq k-1$.

\smallskip

\item $\phi(0\bu) = (10^{k})^\infty$ for any sequence $\bu$ having a prefix of the form 
$0^k1v$ where $k\geq 1$, $|v| = k$, $|v|_0 \leq k-1$. 
\end{enumerate}

\medskip

\item By part (iii)(a) of Theorem~\ref{T:w-from-u}, 
$\phi(0\bu) = (10100100)^\infty = (1w_{3,8}0)^\infty$ for any sequence $\bu$ 
beginning with
\[
Pal(010)011 = w_{3,8}011 = \mathbf{010010}011.
\]

\item Let $\bu$ be the (non-characteristic) Sturmian sequence
\[
1\bbf = \mathbf{101}00101001001010010100100101\cdots
\]
where $\bbf$ is the (binary) Fibonacci sequence (see Examples~\ref{Ex:Fibonacci} and 
\ref{Ex:Sturmian}). Then the longest central prefix of $\bu$ is $w_{3,5} = 101 = Pal(10)$ 
and $\bu$ begins with $w_{3,5}00$. Therefore, by part (iii)(b) of Theorem~\ref{T:w-from-u}, 
we have $\phi(0\bu) = \phi(01\bbf) = (10)^\infty = (1w_{1,2}0)^\infty$.

\medskip

\item By part (iii)(c) of Theorem~\ref{T:w-from-u}, 
$\phi(0\bu) = (10100)^\infty = (1w_{2,5}0)^\infty$ for any sequence $\bu$ 
beginning with
\[
Pal(010)101 = w_{3,8}101 = \mathbf{010010}101.
\]

\item By parts (iii)(b) and (iii)(c) of Theorem~\ref{T:w-from-u}, 
$\phi(\bx) = (10)^\infty = (1w_{1,2}0)^\infty$ for any sequence $\bx$ beginning 
with $0011$ or $0100$. In particular, $\phi(0\bt) = (10)^\infty$ for the {\em Thue-Morse 
sequence} $\bt$, which is the fixed point beginning with $0$ of the morphism $0 \mapsto 01$, 
$1 \mapsto 10$:
\[
\bt = \mathbf{0}110100110010110\cdots
\]
Also note that $\phi(\bt) = (110)^\infty = (1w_{2,3}0)^\infty$.

\medskip

\item Recall that the central word $w_{p,q}$ of slope $p/q \in (0,1)$ (where $\gcd(p,q) = 1$) 
has length $q - 2$ and contains $p - 1$ occurrences of $1$ (and $q - p - 1$ occurrences of $0$). 
We observe that if $p > 2q$ (i.e., if $w_{p,q}$ contains more $1$'s than $0$'s), then 
$w_{p,q}$ begins with $1$; otherwise, if $p < 2q$, then $w_{p,q}$ begins with $0$. Hence, we 
deduce the following general facts from part (iii)(a) of Theorem~\ref{T:w-from-u}.

\smallskip

\begin{enumerate}
\item If $p > 2q$, then $\phi(0\bu) = (1w_{p,q}0)^\infty$ for any sequence $\bu$ 
beginning with $w_{p,q}100$.

\smallskip

\item If $p < 2q$, then $\phi(0\bu) = (1w_{p,q}0)^\infty$ for any sequence $\bu$ 
beginning with $w_{p,q}011$.
\end{enumerate}
\end{enumerate}
\end{example}

\begin{remark} \label{R:find-v}
To determine the longest central prefix of a sequence $\bu \in \{0,1\}^\NN$ (which is neither 
a characteristic Sturmian sequence nor a characteristic periodic balanced sequence), possibly 
the easiest way is to check each palindromic prefix of $\bu$ (in order of increasing length) 
to see if it is equal to $Pal(u)$ for some $u \in \{0,1\}^*$, until there are no more 
palindromic prefixes or until one reaches a palindromic prefix that is not in the image of $Pal$.
\end{remark}

\begin{note} Theorem~\ref{T:w-from-u} also provides a method for computing $F(x)$ in 
Case~(v) of Theorem~\ref{main}. For example, $F(\frac{1}{4}) = \frac{2}{3}$ since the 
base 2 expansion of $\frac{1}{4}$ is $01000 \cdots$  (or $00111 \cdots$) and we have 
$\phi(01000 \cdots) = (10)^\infty = \phi(00111\cdots)$ where $(10)^\infty$ is the 
base 2 expansion of $2/3$. (See part (1) of Example~\ref{Ex:computation} above.)
\end{note}

\section{Larger bases}
What precedes is essentially using base $2$ expansions. One may ask what happens with 
base $b$ expansions, where $b \geq 3$, or what can be said about the intervals 
containing all $\{\xi b^n\}$ for some $\xi$. The result of Bugeaud and Dubickas in 
\cite{yBaD05frac} recalled at the beginning implies that Sturmian sequences (with 
values on an alphabet $\{k, k+1\}$ for some $k \in \{0, 1, \ldots, b-2\}$) will 
again play a fundamental r\^ole.

\section{Acknowledgments}
The second author was partially supported by the Icelandic Research Fund (grant
no. 09003801/1). Both authors thanks the referee for a careful reading and 
useful remarks.

\end{document}